\documentclass{amsart}
\usepackage{amsmath}
\usepackage{amsthm}
\usepackage{amssymb,amscd}
\numberwithin{equation}{section}
\theoremstyle{definition}
\newtheorem{dfn}{Definition}[section]
\newtheorem{example}[dfn]{Example}
\newtheorem{rem}[dfn]{Remark}
\theoremstyle{plain}
\newtheorem{thm}[dfn]{Theorem}
\newtheorem{prop}[dfn]{Proposition}
\newtheorem{cor}[dfn]{Corollary}
\newtheorem{lem}[dfn]{Lemma}
\newtheorem{mainlem}[dfn]{Main Lemma}
\title{Hopf structures on ambiskew polynomial rings}
\author{Jonas T. Hartwig}
\address{Department of Mathematical Sciences, Division of Mathematics,
Chalmers University of Technology and
University of Gothenburg, Eklandagatan 86,
S-412 96 Gothenburg, Sweden}
\email{jonas.hartwig@math.chalmers.se}
\subjclass[2000]{Primary 16S36; Secondary 16W30, 17B37}
\def\car{\mathsf{h}}
\def\a1{\xi}
\def\sl2{\mathfrak{sl_2}}
\def\osp{\mathfrak{osp}}

\def\al{\alpha}
\def\be{\beta}
\def\ga{\gamma}
\def\De{\Delta}
\def\ep{\varepsilon}
\def\ze{\zeta}
\def\la{\lambda}
\def\si{\sigma}

\def\ph{\varphi}
\def\fn{\mathfrak{n}}
\def\fm{\mathfrak{m}}

\def\Z{\mathbb{Z}}
\def\C{\mathbb{C}}
\def\K{\mathbb{K}}
\def\un{\underline}
\DeclareMathOperator{\supp}{Supp} 
 
\DeclareMathOperator{\id}{Id} \DeclareMathOperator{\Max}{Max}
\DeclareMathOperator{\Rad}{Rad}
\begin{document}
\begin{abstract}
We derive necessary and sufficient conditions for an
ambiskew polynomial ring to have a Hopf algebra structure
of a certain type. This construction generalizes many known
Hopf algebras, for example $U(\sl2)$, $U_q(\sl2)$ and
the enveloping algebra of the 3-dimensional Heisenberg Lie algebra.
In a torsion-free case
we describe the finite-dimensional simple modules,
in particular their dimensions and prove a Clebsch-Gordan
decomposition theorem for the tensor product of two simple modules.
We construct a Casimir type operator and
prove that any finite-dimensional weight module is semisimple.
\end{abstract}
\maketitle
\section{Introduction}
In \cite{BDE}, the authors define a four parameter deformation of
the Heisenberg (oscillator) Lie algebra $\mathcal{W}_{\al,\be}^\ga(q)$
and study its representations. Moreover by requiring this algebra to
be invariant under $q\to q^{-1}$, they define a Hopf algebra structure
on $\mathcal{W}_{\al,\be}^\ga(q)$ generalizing several previous
results.

The quantum group $U_q(\mathfrak{sl}_2(\mathbb{C}))$ has by
definition the structure of a Hopf algebra. In
\cite{DQS}, an extension of this quantum group to an associative
algebra denoted by $U_q(f(H,K))$ (where $f$ is a Laurent
polynomial in two variables) is defined and finite-dimensional
representations are studied. The authors show that under certain
conditions on $f$, a Hopf algebra structure can be introduced.
Among these Hopf algebras is for example
the Drinfeld double $\mathcal{D}(\sl2)$.

All of the mentioned algebras fall (after suitable mathematical
formalization in the case of $\mathcal{W}_{\al,\be}^\ga(q)$)
into the class of
so called ambiskew polynomial rings (see Section \ref{sec:prel} for
the definition).
Motivated by these examples of similar classes of algebras, all
of which can be equipped with Hopf algebra structures, we consider
a certain type of Hopf structures on a class of ambiskew polynomial rings.

In Section \ref{sec:prel}, we recall some definitions and fix notation.
We present the conditions for a certain
Hopf structure on an ambiskew polynomial ring in Section \ref{sec:hopf},
while Section \ref{sec:examples} is devoted to examples.
In Section \ref{sec:function} we introduce some convenient
notation and state some useful formulas for
viewing $R$ as an algebra of functions on its set of maximal ideals.
Finite-dimensional simple modules are
studied in Section \ref{sec:simplemodules}.
Those have already been classified in \cite{Jor2}, but we focus on
describing the dimensions in terms of the highest weights. The
main result is stated in Theorem \ref{thm:mainresult_ny}.
The classical Clebsch-Gordan theorem for $U(\sl2)$
is generalized in Section \ref{sec:CG} to the present
more general setting, using the results of the previous section.
Finally, in Section \ref{sec:semisimplicity}
we first construct a kind of Casimir operator and
prove that it can be used to distinguish non-isomorphic
simple modules. This is then used to prove that
any weight module is semisimple.

\section{Preliminaries}\label{sec:prel}
Throughout, $\K$ will be an algebraically closed field of
characteristic zero. All algebras are associative and unital $\K$-algebras.

By a \emph{Hopf structure} on an algebra $A$ we mean a triple $(\De,\ep,S)$
where the \emph{coproduct} $\De:A\to A\otimes A$ is a homomorphism,
($A\otimes A$ is given the tensor product algebra structure)
the \emph{counit} $\ep: A\to \K$ is a homomorphism, and
the \emph{antipode} $S:A\to A$ is an anti-homomorphism such that
\begin{align}
  \label{eq:HopfAxiomCoassoc}
(\id\otimes\De)(\De(x))&=(\De\otimes\id)(\De(x)),
&\text{(Coassociativity)}\\
  \label{eq:HopfAxiomCounit}
(\id\otimes\ep)(\De(x))=&\;x=(\ep\otimes \id)(\De(x)),
&\text{(Counit axiom)}\\
m\Big((S\otimes\id)(\De(x))\Big)=&\;\ep(x)=m\Big((\id\otimes S)(\De(x))\Big),
&\text{(Antipode axiom)}
\end{align}
for all $x\in A$. Here $m:A\otimes A\to A$ denotes the multiplication
map of $A$. A \emph{Hopf algebra} is an algebra equipped
with a Hopf structure. An element $x\in A$ of a Hopf algebra $A$
is called \emph{grouplike} if $\De(x)=x\otimes x$ and \emph{primitive}
if $\De(x)=x\otimes 1+1\otimes x$. In the former case it follows
from the axioms that $\ep(x)=1$, $x$ is invertible and $S(x)=x^{-1}$
while in the latter $\ep(x)=0$ and $S(x)=-x$.

If $V_i$ ($i=1,2$) are two modules over a
Hopf algebra $H$, then $V_1\otimes V_2$ becomes an $H$-module
in the following way
\begin{equation}
  \label{eq:tensordef}
a(v_1\otimes v_2)=\sum_i (a_i'v_1)\otimes (a_i''v_2)
\end{equation}
for $v_i\in V_i$ ($i=1,2$) if $a\in H$ with $\De(a)=\sum_i a_i'\otimes a_i''$.
From \eqref{eq:HopfAxiomCoassoc} it follows that if $V_i$ ($i=1,2,3$)
are modules over $H$ then the natural vector space isomorphism
$V_1\otimes(V_2\otimes V_3)\simeq (V_1\otimes V_2)\otimes V_3$
is an isomorphism of $H$-modules. From \eqref{eq:HopfAxiomCounit}
follows that the one-dimensional module $\K_\ep$ associated to the
representation $\ep$ of $H$ is a tensor unit, i.e. $\K_\ep\otimes V
\simeq V\simeq V\otimes \K_\ep$ as $H$-modules for any $H$-module $V$.

Let $R$ be a finitely generated commutative algebra over $\K$. Let
$\si$ be an automorphism of $R$, $\car\in R$ and
$\a1\in\K\backslash\{0\}$.
Then we define the algebra $A=A(R,\si,\car,\a1)$ as the
associative $\K$-algebra formed by adjoining to $R$ two symbols
$X_+$, $X_-$ subject to the relations
\begin{equation}\label{eq:caserel1}
X_\pm a = \si^{\pm 1}(a) X_\pm\text{ for }a\in R,
\end{equation}
\begin{equation}\label{eq:caserel2}
X_+X_- = \car+\a1 X_-X_+.
\end{equation}
This algebra is called an \emph{ambiskew polynomial ring}. Its
structure and representations were studied by Jordan \cite{Jor1}
(see also references therein).

We recall the definition of a \emph{generalized Weyl algebra (GWA)}
(see \cite{Bav} and references therein). If $B$ is a ring, $\si$ an
automorphism of $B$, and $t\in B$ a central element, then the
generalized Weyl algebra $B(\si,t)$ is the ring extension of $B$
generated by two elements $x_+$, $x_-$ subject to the relations
\begin{equation}
\begin{split}
x_\pm a=\si^{\pm 1}(a)x_\pm,\quad\text{for }a\in B,\\
x_-x_+ = t, \quad\text{and}\quad x_+x_- = \si(t).
\end{split}
\end{equation}
The relation between these two constructions is the following.
Let $A=A(R,\si,\car,\a1)$ be an ambiskew polynomial ring.
Denote by $R[t]$ be the polynomial ring in one variable $t$ with
coefficients in $R$ and let us
extend the automorphism $\si$ of $R$ to a $\K$-algebra
automorphism of $R[t]$ satisfying
\begin{equation}
  \label{eq:extaut}
\si(t)=\car+\a1 t.
\end{equation}
Then $A$ is isomorphic to the GWA $R[t](\si,t)$.

%
\section{The Hopf structure} \label{sec:hopf}
%
Let $A=A(R,\si,\car,\a1)$ be a skew polynomial ring and
assume that $R$ has been equipped with a Hopf structure.
In this section we will extend the Hopf structure on $R$
to $A$. We make the following anzats, guided by \cite{BDE} and
\cite{DQS}:
\begin{align}
\De(X_\pm)&=X_\pm\otimes r_\pm + l_\pm\otimes X_\pm,\label{eq:ansatzDe}\\
\ep(X_\pm)&=0,\label{eq:ansatzep}\\
S(X_\pm)&=s_\pm X_\pm.\label{eq:ansatzS}
\end{align}
The elements $r_\pm,l_\pm$ and $s_\pm$ will be assumed to belong
to $R$.
\begin{thm}\label{thm:Hopf}
Formulas (\ref{eq:ansatzDe})-(\ref{eq:ansatzS}) define a Hopf
algebra structure on $A$ which extends that of $R$ iff
\begin{subequations}\label{eq:thmHopfFIRST}
\begin{align}
\label{eq:thmHopfsiDe}
(\si\otimes\id)\circ\De |_{R} &=\De\circ\si|_{R}=(\id\otimes\si)\circ\De|_{R},\\
\label{eq:TH_S1} S\circ\si |_{R} &= \si^{-1}\circ S|_{R},
\end{align}
\end{subequations}
\begin{subequations}
\begin{align}
\label{eq:thmHopfDe1}
\De(\car)&=\car\otimes r_+r_-+l_+l_-\otimes \car,\\
\ep(\car)&=0, \label{eq:TH_epcartan}\\
\label{eq:TH_S3} S(\car) &= -(l_+l_-r_+r_-)^{-1}\car,
\end{align}
\end{subequations}
\begin{subequations}
\begin{align}
\label{eq:TH_grouplike}
r_\pm \text{ and } l_\pm \text{ are grouplike, i.e. }&\text{$\De(x)=x\otimes x$ for $x\in\{r_\pm,l_\pm\}$,}\\
\label{eq:thmHopfDe23} \si(l_\pm)\otimes\si(r_\mp)&=\a1 l_\pm\otimes r_\mp,
\end{align}
\end{subequations}
\begin{equation} \label{eq:TH_s}
 (s_\pm)^{-1}=-l_\pm\si^{\pm 1}(r_\pm).
\end{equation}
%
\end{thm}
\begin{proof}
From \eqref{eq:caserel1}-\eqref{eq:caserel2} we see that
$\ep$ extends to a homomorphism $A\to\K$ satisfying \eqref{eq:ansatzep} if and only if
\eqref{eq:TH_epcartan} holds.
Assume for a moment that $\De$ extends to a homomorphism $A\to A\otimes A$.
From (\ref{eq:ansatzDe})-(\ref{eq:ansatzep}) it
follows that $\ep$ is a counit iff
\begin{equation} \label{eq:epr=1}
\ep(r_+)=\ep(r_-)=\ep(l_+)=\ep(l_-)=1.
\end{equation}
$\De$ is coassociative iff (dropping the $\pm$)
\[(\id\otimes\De)(\De(X))=(\De\otimes\id)(\De(X))\]
which is equivalent to
\[X\otimes \De(r)+l\otimes X\otimes r + l\otimes l\otimes X=
X\otimes r\otimes r+l\otimes X\otimes r + \De(l)\otimes  X, \]
or
\begin{equation}
  \label{eq:flodmyra}
X\otimes (\De(r)-r\otimes r) = (\De(l)-l\otimes l)\otimes X.
\end{equation}
From \eqref{eq:caserel1}-\eqref{eq:caserel2} follows that
$A$ has a $\Z$-gradation defined by requiring
that $\deg r=0$ for $r\in R$, $\deg X_\pm=\pm 1$.
This also induces a $\Z^2$-gradation on $A\otimes A$ in a natural
way.
The left and right hand sides of equation \eqref{eq:flodmyra} are
homogenous of different
$\Z^2$-degrees, namely $(\pm 1,0)$ and $(0,\pm 1)$ respectively.
Hence, since homogenous elements of different
degrees must be linearly independent,
\eqref{eq:flodmyra} is equivalent to both sides being zero
which holds iff $r_\pm$ and $l_\pm$ are grouplike.

$\De$ respects \eqref{eq:caserel1} iff (again dropping $\pm$)
\begin{align*}
\De(X) \De(a)&=\De(\si(a))\De(X),\\
(X\otimes r+l\otimes X)\De(a)&= \De(\si(a))(X\otimes r+l\otimes
X),
\end{align*}
\begin{equation*}
(\si\otimes 1)(\De(a))\cdot(X\otimes r)+
(1\otimes\si)(\De(a))\cdot(l\otimes X)= \De(\si(a))(X\otimes
r+l\otimes X),
\end{equation*}
\begin{equation*}
\big((\si\otimes 1)(\De(a))- \De(\si(a))\big) \cdot(X\otimes r)
+\big((1\otimes\si)(\De(a))-\De(\si(a))\big) \cdot(l\otimes X)=0.
\end{equation*}
As before the two terms in the last equation
have different $\Z^2$-degrees and therefore
must be zero. So $\De$ respects \eqref{eq:caserel1} iff
\eqref{eq:thmHopfsiDe} holds.

It is straightforward to check that $\De$ respects \eqref{eq:caserel2} iff
\begin{multline}\label{eq:myra3}
\car\otimes r_+r_-+l_+l_-\otimes \car-\De(\car)+\\
+\Big(l_+\otimes \si(r_-)-\a1 \si^{-1}(l_+)\otimes r_-\Big)X_-\otimes X_++\\
+\Big(\si(l_-)\otimes-\a1 l_-\otimes\si^{-1}(r_+)\Big)X_+\otimes X_-=0.
\end{multline}
Again these three terms have different degrees so each of them
must be zero. Hence \eqref{eq:thmHopfDe1} holds.
Multiply the second term by $X_+\otimes X_-$ from the right:
\[\big(l_+\otimes \si(r_-)-\a1 \si^{-1}(l_+)\otimes r_-\big)t\otimes \si(t)=0.\]
Here we use the extension \eqref{eq:extaut} of
$\si$ to $R[t]$ where $t=X_-X_+$.
If we apply $e_1\otimes e_1'$ to this equation, where $e_r$
($e_r'$) for $r\in R$ is the evaluation homomorphism $R[t]\to R$
which maps $t$ ($\si(t)$) to $r$, we get
\[l_+\otimes \si(r_-)=\a1 \si^{-1}(l_+)\otimes r_-.\]
Applying $\si\otimes 1$ to this we obtain one of the relations in \eqref{eq:thmHopfDe23}.
Similarly the vanishing of the third term in \eqref{eq:myra3}
implies the other.

Assuming that $S$ is an anti-homomorphism $A\to A$ satisfying
\eqref{eq:ansatzS}, we obtain that $S$ is a  an antipode on $A$ iff
\[S(X_\pm)r_\pm+S(l_\pm)X_\pm=0=X_\pm S(r_\pm)+l_\pm S(X_\pm),\]
which is equivalent to \eqref{eq:TH_s}, using that $r_\pm$ and $l_\pm$
are grouplike.
And $S$ extends to a well-defined anti-homomorphism $A\to A$ iff
\begin{align}
  \label{eq:THP_S1}
  S(a)S(X_\pm)&=S(X_\pm)S(\si^{\pm 1}(a)),\quad \text{for $a\in R$,}\\
  \label{eq:THP_S2}
  S(X_-)S(X_+)&=S(\car)+\a1 S(X_+)S(X_-).
\end{align}
Using \eqref{eq:TH_s} and that $r_\pm,l_\pm$ are invertible,
\eqref{eq:THP_S1} holds iff \eqref{eq:TH_S1} holds.
And \eqref{eq:THP_S2} holds iff
\begin{align*}
0&=s_-X_-s_+X_+-S(\car)-\a1 s_+X_+s_-X_-=\\
&=s_-\si^{-1}(s_+)X_-X_+-S(\car)-s_+\a1 \si(s_-)X_+X_-=\\
&=-S(\car)-s_+\si(s_-)\a1 \car+\\
&\quad+\big(s_-\si^{-1}(s_+)-s_+\si(s_-)\a1^2\big)t.
\end{align*}
Applying $e_0$ and $e_1$ we obtain
\begin{align*}
S(\car)&=-\a1 s_+\si(s_-) \car,\\
s_-\si^{-1}(s_+)&=\a1^2 s_+\si(s_-).
\end{align*}
Substituting \eqref{eq:TH_s} in these equations and using
\eqref{eq:thmHopfDe23}, the first is equivalent to \eqref{eq:TH_S3},
while the other already holds.
\end{proof}

%
%
\section{Examples} \label{sec:examples}
Many Hopf algebras known in the literature can be viewed as one
defined in the previous section.
\subsection{Heisenberg algebra}\label{sec:U(h_3)}
Let $R=\C[c]$ with $c$ primitive, and $\si(c)=c$. Choose $\car=c$,
$\a1=r_+=r_-=l_+=l_-=1$. Then $A$ is the universal enveloping algebra
$U(\mathfrak{h}_3)$ of the
three-dimensional Heisenberg Lie algebra.

\subsection{$U(\sl2)$ and its quantizations}
\subsubsection{$U(\sl2)$}\label{sec:Usl2}
Let $R=\C[H]$ with Hopf algebra structure $\De(H)=H\otimes
1+1\otimes H$, $\ep(H)=0$, $S(H)=-H$. Define $\si(H)=H-1$. Choose
$\car=H$, $\a1=r_+=r_-=l_+=l_-=1$. Then $A\simeq U(\sl2)$ as Hopf
algebras.

\subsubsection{$U_q(\sl2)$}\label{sec:Uqsl2}
Let $R=\C[K,K^{-1}]$ with Hopf structure defined by requiring that
$K$ is grouplike. Define $\si(K)=q^{-2}K$, where $q\in\C,
q^2\neq 1$, and choose $\car=\frac{K-K^{-1}}{q-q^{-1}}$,
$\a1=r_-=l_+=1$ and $r_+=K$, $l_-=K^{-1}$. Then the equations in
Theorem \ref{thm:Hopf} are satisfied giving a Hopf algebra
$A$ which is isomorphic to $U_q(\sl2)$.

\subsubsection{$\breve{U}_q(\sl2)$}
For the definition of this algebra,
see for example \cite{KliSch}. Let $q\in\C,
q^4\neq 1$. Let $R=\C[K,K^{-1}]$ with $K$ grouplike. Define
$\si(K)=q^{-1}K$, $\car=\frac{K^2-K^{-2}}{q-q^{-1}}$, $\a1=1$,
$r_+=r_-=K$, $l_+=l_-=K^{-1}$. Then $A=A(R,\si,\car,\a1)$ is a Hopf algebra
isomorphic to $\breve{U}_q(\sl2)$.

\subsection{$U_q(f(H,K))$} \label{sec:dqs}
Let $R=\C[H,H^{-1},K,K^{-1}]$, $\si(H)=q^2H$, $\si(K)=q^{-2}K$.
Let $\al\in\K$ and $M,p,r,s,t,p',r',s',t'\in\Z$ such that
$M=m-n=m'-n'=p+t-r-s$, $s-t=s'-t'$ and $p-r=p'-r'$. Set $\car=\al
(K^mH^n-K^{-m'}H^{-n'})$, $\a1=1$, $r_+=K^pH^r$, $l_+=K^sH^t$,
$r_-=K^{-s'}H^{-t'}$, $l_-=K^{-p'}H^{-r'}$. Then $A$ is the Hopf
algebra described in \cite{DQS}, Theorem 3.3.

\subsection{Down-up algebras} \label{sec:downup}
The down-up algebra $A(\al,\be,\ga)$ where $\al,\be,\ga\in \C$,
was defined in \cite{BenRob} and studied by many authors, see for example
\cite{BenWit}, \cite{CarMus}, \cite{Jor1}, \cite{KirMus},  and references
therein. It is the algebra generated by $u,d$ and relations
\begin{align*}
  ddu&=\al dud+\be udd+\ga d,\\
  duu&=\al udu+\be uud+\ga u.
\end{align*}
In \cite{Jor1} it is proved that if $\si$ is allowed to be
any endomorphism, not necessarily invertible, then
any down-up algebra is an ambiskew polynomial ring.
Here we consider the down-up algebra $B=A(0,1,1)$.
Thus $B$ is the $\C$-algebra with generators $u,d$ and relations
\begin{equation}
  \label{eq:A011rels}
d^2u=ud^2+d,\quad du^2=u^2d+u.
\end{equation}
Let $R=\C[\car]$, $\si(\car)=\car+1$ and $\xi=-1$.
Then $B$ is isomorphic to the ambiskew polynomial
ring $A(R,\si,\car,\xi)$ via $d\mapsto X_+$ and $u\mapsto X_-$.

One can show that $B$ is isomorphic to the enveloping algebra of
the Lie super algebra $\osp(1,2)$ and hence has a \emph{graded}
Hopf structure. A question was raised in \cite{KirMus} whether there
exists a Hopf structure on $B$. We do not answer this question
here but we show the existence of a
Hopf structure on a larger algebra $B_q$ giving us a formula
for the tensor product of weight (in particular finite-dimensional)
modules over $B$.

Let $q\in\C^*$ and fix a value of $\log q$. By $q^a$
we always mean $e^{a\log q}$.
Let $B_q$ be the ambiskew polynomial ring
$B_q=A(R,\si,\car,\xi)$ where $R=\C[\car,w,w^{-1}]$,
$\si(\car)=\car+1$, $\si(w)=qw$, and $\xi=-1$.

\begin{thm}
For any $\rho,\la\in\Z$ such that
$q^{\rho-\la}=-1$ and $q^{2\rho}=1$, the algebra $B_q$ has a Hopf algebra
structure given by
\begin{align*}
  \De(X_\pm)&=X_\pm\otimes w^{\pm\rho}+w^{\pm\la}\otimes X_{\pm},\quad
  \ep(X_\pm)=0,\\
  S(X_\pm)&=-w^{\mp\la}X_\pm w^{\mp\rho}=-q^{\rho}w^{\mp(\rho+\la)}X_\pm,
\intertext{and}
\De(w)&=w\otimes w, \quad \ep(w)=1,\quad S(w)=w^{-1},\\
\De(\car)&=\car\otimes 1+1\otimes \car,\quad \ep(\car)=0,\quad S(\car)=-\car.
\end{align*}
\end{thm}
\begin{proof}
The subalgebra $\C[\car,w,w^{-1}]$ of $B_q$ has a unique Hopf
structure given by the maps above.
We must verify \eqref{eq:thmHopfFIRST}-\eqref{eq:TH_s}
with $\car = v$, $\xi=-1$, $r_\pm = w^{\pm\rho}$, $l_\pm=w^{\pm\la}$,
and $s_\pm=-q^{\rho}w^{\mp(\rho+\la)}$. This is straightforward.
\end{proof}

This gives us a tensor structure on the category of
modules over $B_q$. Next aim is to show how using
the Hopf structure on $B_q$ one can define a tensor structure on the
category of weight modules over $B$.

In general, if $C$ is
a commutative subalgebra of an algebra $A$, we say that
an $A$-module $V$ is a \emph{weight module} with respect to $C$ if
\[V=\oplus_{\fm\in \Max(C)}V_\fm,\qquad V_\fm=\{v\in V|\fm v=0\},\]
where $\Max(C)$ denotes the set of all maximal ideals of $C$.
When $C$ is finitely generated this is equivalent to
$V$ having a basis in which each $c\in C$ acts diagonally.

By weight modules over $B$ ($B_q$) we mean weight modules with respect to
the subalgebra $\C[\car]$ ($\C[\car,w,w^{-1}]$). We need a simple lemma.
\begin{lem}\label{lem:Bfdweight}
Any finite-dimensional module $V$ over $B$ is a weight module.
\end{lem}
\begin{proof}
By Proposition 5.3 in \cite{Jor1}, any finite-dimensional $B$-module
is semisimple. Since direct sums of weight modules are weight modules we
can assume that $V$ is simple. Since $V$ is finite-dimensional,
the commutative subalgebra $\C[\car]$ has a common eigenvector $v\neq 0$,
i.e. $\fm v=0$ for some maximal ideal $\fm$ of $\C[\car]$.
Acting on this weight vector by $X_\pm$ produces another weight vector:
$\si^{\pm 1}(\fm)X_\pm v=X_\pm\fm v=0$. Since
$B$ is generated by $\C[\car]$ and $X_\pm$, any vector
in the $B$-submodule of $V$
generated by $v$ is a sum of weight vectors. But $V$ was
simple so $V=\oplus_\fm V_\fm$.
\end{proof}

Let $\mathcal{W}(B)$ denote the category of weight $B$-modules
and similarly for $B_q$.
\begin{thm}
The category of weight modules over $B$ can be embedded into the
category of weight modules over $B_q$, i.e. there exist functors
\[\mathcal{W}(B)\overset{\mathcal{E}}{\longrightarrow}
\mathcal{W}(B_q)
\overset{\mathcal{R}}{\longrightarrow} \mathcal{W}(B) \]
whose composition is the identity functor.
In particular, the category of finite-dimensional $B$-modules
can be embedded in $\mathcal{W}(B_q)$.
\end{thm}
\begin{proof}
$\mathcal{R}$ is given by restriction. It takes weight modules
to weight modules. Next we define $\mathcal{E}$.
Let $V$ be a weight module over $B$ and define
\begin{equation}
  \label{eq:waction}
wv=q^{\al}v\qquad\text{ for $v\in V_{(\car-\al)}$ and $\al\in\C$.}
\end{equation}
It is immediate that $w$ commutes with $\car$. Let $v\in V_{(\car-\al)}$
be arbitrary. Then
\[X_+wv=X_+q^{\al}v=q^{\al}X_+v.\]
On the other hand, since $\car X_+v = X_+(\car-1)v=(\al-1)X_+v$
which shows that $X_+v\in V_{(\car-(\al-1))}$, we have
\[qwX_+v=qq^{\al-1}X_+v=q^{\al}X_+v.\]
Thus $X_+w=qwX_+$.
Similarly $X_-w=q^{-1}wX_-$ on $V$. Thus $V$ becomes
a module over $B_q$. That $V$ is a weight module with respect
to $\C[\car,w,w^{-1}]$ is clear. We define $\mathcal{E}(V)$ to be the
same space $V$ with additional action \eqref{eq:waction}.
If $\ph:V\to W$ is a morphism
of weight $B$-modules then $\ph(wv)=w\ph(v)$ for weight vectors $v$,
since $\ph(V_\fm)\subseteq W_\fm$
for any maximal ideal $\fm$ of $\C[\car]$. But then
$\ph(wv)=w\ph(v)$ for all $v\in V$ since $V$ is a weight module.
Thus $\ph$ is automatically
a morphism of $B_q$-modules and we set $\mathcal{E}(\ph)=\ph$.
It is clear that the composition of the functors is the
identity on objects and morphisms.
\end{proof}

Note that
\begin{equation}
  \label{eq:image}
  \mathcal{E}\big(\mathcal{W}(B)\big)=\big\{V\in \mathcal{W}(B_q) \; |\;
\supp(V)\subseteq \{\fm=(\car-\al,w-q^\al) \; |\; \al\in\C\}\big\}.
\end{equation}

It is not difficult to see that
\[\mathcal{E}(V_1)\otimes\mathcal{E}(V_2)\in
\mathcal{E}\big(\mathcal{W}(B)\big)\]
and hence there is a unique $V_3\in \mathcal{W}(B)$ such that
\[\mathcal{E}(V_1)\otimes\mathcal{E}(V_2)=\mathcal{E}(V_3).\]
Thus we can define
\[V_1\otimes V_2:=V_3\]
and this will make $\mathcal{W}(B)$ into a tensor category.

\subsection{Non Hopf ambiskew polynomial rings}
There are many examples of ambiskew polynomial rings
which do not have any
Hopf structure. One example is the Weyl algebra $W=\langle a,b| ab-ba=1\rangle$
which can have no counit $\ep$. Indeed, a counit is in particular a
homomorphism $\ep:W\to\C$
so we would have $1=\ep(1)=\ep(a)\ep(b)-\ep(b)\ep(a)=0$.
Moreover all down-up algebras are ambiskew polynomial rings
(see \cite{Jor1}) and \cite{KirMus} contains necessary conditions for the existence
of a Hopf structure on a down-up algebra in terms of the
parameters $\al,\be,\ga$. More precisely, they show that if $A=A(\al,\be,\ga)$
is a Noetherian down-up algebra that is a Hopf algebra, then $\al+\be=1$.
Moreover if $\ga=0$, then $(\al,\be)=(2,-1)$ and as algebras, $A$ is
isomorphic to the universal enveloping algebra of the three-dimensional
Heisenberg Lie algebra, while if $\ga\neq 0$, then $-\be$ is not an $n$th root
of unity for $n\ge 3$. It would be of interest to generalize
such a result to a more general class of ambiskew polynomial rings
and also to other GWAs.

%
%
\section{$R$ as functions on a group} \label{sec:function}
From now on we assume that $A=A(R,\si,\car,\xi)$ is an algebra of
the form defined in Section \ref{sec:hopf} and that conditions
\eqref{eq:thmHopfFIRST}-\eqref{eq:TH_s} hold so that $A$ becomes a
Hopf algebra with $R$ as a Hopf subalgebra. Let $G$ denote the set
of all maximal ideals in $R$. Since $\K$ is algebraically closed
and $R$ is finitely generated, the inclusion map $i_\fm:\K\to
R/\fm$ is onto for any $\fm\in G$ and we let $\ph_\fm:R\to \K$
denote the composition of the projection $R\to R/\fm$ and
$i_\fm^{-1}$. Thus $\ph_\fm(a)$ is the unique element of $\K$ such that
$a-\ph_\fm(a)\in \fm$. We define the \emph{weight sum} of
$\fm,\fn\in G$ to be
\[\fm+\fn := \ker (m\circ (\ph_\fm\otimes\ph_\fn)\circ\De|_{R}).\]
This is the kernel of a $\K$-algebra homomorphisms $R\to \K$,
hence $\fm+\fn\in G$. We will never use the usual addition of
ideals so $+$ should not cause any confusion. Using that $\De$ is
coassociative, $\ep$ is a counit and $S$ is an antipode, one
easily deduces that $+$ is associative, that $\un{0}:=\ker \ep$ is
a unit element and $S(\fm)$ is the inverse of $\fm$. Thus $G$ is a
group under $+$. If $R$ is cocommutative, $G$ is abelian.

\begin{example} Let $R=\C[H]$. Then $G=\{(H-\al)\,|\,\al\in\C\}$. Give
$R$ the Hopf structure $\De(H)=H\otimes 1+1\otimes H$, $\ep(H)=0$ and
$S(H)=-H$. Then the operation $+$ will be
\[(H-\al)+(H-\be)=(H-(\al+\be)),\]
i.e. the correspondence $\C\ni \al\mapsto (H-\al)\in G$ is an additive group
isomorphism.

 If $R=\C[K,K^{-1}]$ then $G=\{(K-\al)\,|\,\al\in\C^*\}$. With
the Hopf structure $\De(K)=K\otimes K$, $\ep(K)=1$ and $S(K)=K^{-1}$,
the operation $+$ will be
\[(K-\al)+(K-\be)=(K-\al\be)\]
for $\al,\be\neq 0$. Thus $G\simeq\langle \C^*,\cdot\rangle$.
\end{example}

We will often think of elements from $R$ as $\K$-valued functions
on $G$ and for $x\in R$ and $\fm\in G$ we will use the notation
$x(\fm)$ for $\ph_\fm(x)$. Note however that different elements
$x,y\in R$ can represent the same function. In fact one can check
that the map from $R$ to functions on $G$ is a homomorphism of
$\K$-algebras with kernel equal to the radical
$\Rad(R):=\cap_{\fm\in G} \fm$.

Define a map
\begin{equation}\label{eq:zetadef}
\ze:\Z\to G,\; n\mapsto \un{n}:=\si^n(\un{0}).
\end{equation}
\begin{lem} Let $\fm,\fn\in G$. Then for any $a\in R$,
\begin{align}
\si(a)(\fm)&=a\big(\si^{-1}(\fm)\big),\label{eq:grplem1}\\
a(\fm+\fn)&=m\circ (\ph_\fm\otimes\ph_\fn)\circ\De(a)=\sum_{(a)}a'(\fm)a''(\fn),\label{eq:grplem2}\\
\fm+\un{1}&=\si(\fm) = \un{1}+\fm.\label{eq:grplem3}
\end{align}
Thus $\ze$ is a group homomorphism and
its image is contained in the center of $G$.
\end{lem}
\begin{proof}
Since for any $a\in R$ we have
\[\si(a)(\fm)-a=\si^{-1}\big(\si(a)(\fm)-\si(a)\big)\in\si^{-1}(\fm),\]
\eqref{eq:grplem1} holds. Similarly,
\[a(\fm+\fn)-a\in \fm+\fn\]
so applying the map $m\circ (\ph_\fm\otimes\ph_\fn)\circ\De$
to $a(\fm+\fn)-a$ yields zero. This gives \eqref{eq:grplem2}.
Finally we have for any $a\in\fm$,
\begin{multline*}
\si(a)(\fm+\un{1})=
m\circ (\ph_\fm\otimes\ph_{\un{1}})\circ\De(\si(a))=
m\circ (\ph_\fm\otimes\ph_{\un{1}})\circ(1\otimes \si)\De(a)=\\
=m\circ (\ph_\fm\otimes\ph_{\un{0}})\circ\De(a)=a(\fm+\un{0})=a(\fm)=0.
\end{multline*}
Here we used \eqref{eq:grplem2} in the first and the fourth equality,
\eqref{eq:thmHopfsiDe} in the second and \eqref{eq:grplem1} in the third.
Thus $\si(\fm)\subseteq \fm+\un{1}$ and then equality holds since both
sides are maximal ideals. The proof of the other equality in \eqref{eq:grplem3}
is symmetric.
\end{proof}

\begin{example}  If $R=\C[K,K^{-1}]$ with
 $\De(K)=K\otimes K, \ep(K)=1, S(K)=K^{-1}$ and $\si(K)=q^{-2}K$, then
$\ker\ep=(K-1)$ so
\[\un{n}=\si^n(\un{0})=\si^n((K-1))=(q^{-2n}K-1)=(K-q^{2n}).\]
\end{example}

From \eqref{eq:grplem2} follows that if $x\in R$ is grouplike,
then viewed as a function $G\to \K$ it is a multiplicative
homomorphism. Using \eqref{eq:grplem2} and
\eqref{eq:thmHopfDe1}-\eqref{eq:TH_S3}, the following formulas are
satisfied by $\car$ as a function on $G$.
\begin{equation}\label{eq:hconds}
\begin{split}
\car(\fm+\fn)&=\car(\fm)r(\fn)+l(\fm)\car(\fn),\\
\car(\un{0})&=0,\\
\car(-\fm)&=-r^{-1}l^{-1}\car(\fm),
\end{split}
\end{equation}
where $r=r_+r_-$ and $l=l_+l_-$.



\section{Finite-dimensional simple modules} \label{sec:simplemodules}
In this section we consider finite-dimensional simple modules
over the algebra $A$. The main theorem is Theorem \ref{thm:mainresult_ny}
where we, under the torsion-free assumption \eqref{eq:torfree},
characterize the finite-dimensional simple modules of a given
dimension in terms
of their highest weights. This result will be used in Section \ref{sec:CG}
to prove a Clebsch-Gordan decomposition theorem.

Throughout the rest of the paper we will assume that
\begin{equation}
  \label{eq:torfree}
  \si^n(\fm)\neq\fm\text{ for any $n\in\Z\backslash\{0\}$ and any $\fm\in G$.}
\end{equation}
By \eqref{eq:grplem3}, this condition holds iff
$\un{1}$  has infinite order in $G$.

\subsection{Weight modules, Verma modules and their finite-dimensional
simple quotients}
In this section we define weight modules, Verma modules and derive an equation
for the dimension of its finite-dimensional simple quotients.

Let $V$ be an $A$-module. We call $\fm\in G$
a \emph{weight} of $V$ if $\fm v=0$ for some nonzero $v\in V$.
The \emph{support} of $V$, denoted $\supp(V)$, is the set of weights of $V$.
To a weight $\fm$ we associate its \emph{weight space}
\[ V_\fm=\{v\in V \, |\, \fm v=0\}.\]
Elements of $V_\fm$ are called \emph{weight vectors of weight $\fm$}.
A module $V$ is a \emph{weight module} if $V=\oplus_\fm V_\fm$.
A \emph{highest weight vector} $v\in V$ of weight $\fm$ is a weight vector of
weight $\fm$ such that $X_+v=0$.
A module $V$ is called a \emph{highest weight module} if it is generated by
a highest weight vector. From the defining relations of $A$ it follows
that
\begin{equation}
  \label{eq:weightaction}
X_\pm V_\fm\subseteq V_{\si^{\pm 1}(\fm)}.
\end{equation}
Equation \eqref{eq:weightaction} implies that a highest weight module
is a weight module.

Let $\fm\in G$. The \emph{Verma module} $M(\fm)$ is defined
as the left $A$-module $A/I(\fm)$ where $I(\fm)$ is the
left ideal $AX_++A\fm\subseteq A$.
From relations \eqref{eq:caserel1},\eqref{eq:caserel2}
follows that
\[\{v_n:=X_-^n+I(\fm) \, |\, n\ge 0\}\]
is a basis for $M(\fm)$.
It is clear that $M(\fm)$ is a highest weight module generated by
$v_0$. We also see that
the vectors $v_n\; (n\ge 0)$ are weight vectors of weights $\si^n(\fm)$
respectively.
By \eqref{eq:torfree} we conclude $\dim M(\fm)_\fm=1$.
Therefore the sum of all its proper submodules is proper and equals
the unique maximal submodule $N(\fm)$ of $M(\fm)$. Thus $M(\fm)$
has a unique simple quotient $L(\fm)$.
Since it is easy to see that any highest weight module over $A$
of highest weight $\fm$ is a quotient of $M(\fm)$ we deduce that
$L(\fm)$ is the unique irreducible highest weight module over $A$ with given
highest weight $\fm\in G$.
We set
\[G_f:=\{\fm\in G\; |\; \dim L(\fm)<\infty\}.\]

\begin{prop}\label{prop:anyfd}
Any finite-dimensional simple module
over $A$ is isomorphic to $L(\fm)$ for some $\fm\in G_f$.
\end{prop}
\begin{proof}
Let $V$ be a finite-dimensional simple $A$-module.
Since $\K$ is algebraically closed, $R$ has a common eigenvector
$v\neq 0$, i.e. there exists $\fn\in G$ such that $\fn v=0$. From
\eqref{eq:caserel1} it follows that $\si^n(\fn)(X_+)^nv=0$ for any
$n\ge 0$. By \eqref{eq:torfree}, the set $\{X_+^nv \, |\, n\ge 0\}$ is
a set of weight vectors of different weights.
Since $V$ is finite-dimensional it follows that
$(X_+)^nv=0$ for some $n>0$.
This proves the existence of a highest weight vector of weight
$\fm$ in $V$ for some weight $\fm$. Thus $V=L(\fm)$.
\end{proof}

\begin{cor}\label{cor:weightsupp}
Let $V$ be a finite-dimensional weight module over $A$. Then
\[\supp(V)\subseteq G_f+\un{\Z}=\{\fm+\un{n} \, |\, \fm\in G_f, n\in\Z\}.\]
\end{cor}
\begin{proof} Let $\fm\in\supp(V)$ and let $0\neq v\in V_\fm$.
Then $(X_+)^nv=0$ for some smallest $n>0$. But then $(X_+)^{n-1}v$ is a
highest weight vector so its weight $\si^{n-1}(\fm)=\fm+\un{n-1}$
must belong to $G_f$. Thus $\fm=\fm+\un{n-1}-\un{n-1}\in G_f+\un{\Z}$.
\end{proof}

The following lemma was essentially proved in \cite{Jor2}, Proposition 2.3,
and the general result was mentioned in \cite{Jor1}. We give
a proof for completeness.
\begin{prop} \label{prop:dimLm}
The dimension of $L(\fm)$ is the smallest positive integer
$n$ such that
\[\sum_{k=0}^{n-1}\xi^{n-1-k} \car(\fm-\un{k})=0.\]
\end{prop}
\begin{proof}
Let $e^\fm$ be a highest weight vector in $L(\fm)$. Let $n>0$ be
the smallest positive integer such that $X_-^ne^\fm=0$. Then the
set spanned by the vectors $X_-^je^\fm$, $0\le j<n$, is invariant
under $X_-$, under $R$ using \eqref{eq:caserel1}, and under $X_+$,
using \eqref{eq:caserel2}. Hence it is a nonzero submodule and so
coincides with $L(\fm)$ since the latter is simple. Therefore $n=\dim
L(\fm)$. Let $k>0$. Then $X_-^ke^\fm=0$ implies that
$X_+^kX_-^ke^\fm=0$. Conversely, suppose
$X_+^kX_-^ke^\fm=0$. Then $X_+^{k-1}X_-^ke^\fm$ generates
a proper submodule and thus is zero. Repeating this argument we
obtain $X_-^ke^\fm=0$. Hence $\dim L(\fm)$ is the smallest
positive integer $n$ such that $X_+^nX_-^ne^\fm=0$. Using
induction it is easy to deduce the formulas
\[X_+X_-^n=X_-^{n-1}
\Big(\xi^nX_-X_++\sum_{k=0}^{n-1}\xi^{n-1-k}\si^k(\car)\Big),\]
\begin{equation}\label{eq:xnxn}
X_+^nX_-^n=\prod_{m=1}^n
\Big(\xi^mX_-X_++\sum_{k=0}^{m-1}\xi^{m-1-k}\si^k(\car)\Big).
\end{equation}
Applying both sides of this equality to the vector $e^\fm$ gives
\begin{equation}\label{eq:lemsum}
X_+^nX_-^ne^\fm=\prod_{m=1}^n\sum_{k=0}^{m-1}\xi^{m-1-k}\si^k(h)e^\fm.
\end{equation}
Using that $e^\fm$ is a weight vector of weight $\fm$
and formula \eqref{eq:grplem1} we have
\[\si^k(\car)e^\fm=\si^k(\car)(\fm)e^\fm=\car(\fm-\un{k})e^\fm.\]
Substituting this into \eqref{eq:lemsum} we obtain
\[X_+^nX_-^ne^\fm=\prod_{m=1}^n
\sum_{k=0}^{m-1}\xi^{m-1-k}\car(\fm-\un{k})e^\fm.\]
The smallest positive $n$ such that this is zero must be the
one such that the last factor is zero. The claim is proved.
\end{proof}

\begin{cor}\label{cor:period}
If $\fm,\fm_0\in G$ where $\car(\fm_0)=0$, then
\[\dim L(\fm_0+\fm)=\dim L(\fm)=\dim L(\fm+\fm_0).\]
\end{cor}
\begin{proof} Note that \eqref{eq:hconds} implies
that $\car(\fn+\fm_0)=\car(\fn)r(\fm_0)$ and
$\car(\fm_0+\fn)=l(\fm_0)\car(\fn)$ for any $\fn\in G$,
recall that $r$ and $l$ are invertible and use
Proposition \ref{prop:dimLm}.
\end{proof}
\subsection{Dimension and highest weights}
The goal in this subsection is to prove Theorem \ref{thm:mainresult_ny}
which describes in detail the relationship between
the dimension of a finite-dimensional simple module and its
highest weight.

We begin with a few useful lemmas.
Recall that $r=r_+r_-$ and $l=l_+l_-$. For brevity we set
$r_1=r(\un{1})$ and $l_1=l(\un{1})$. Since $r_\pm,l_\pm$ are grouplike
so are $r$ and $l$ and thus $r_1,l_1$ are nonzero scalars.
\begin{lem} \label{lem:dimabc}
We have\\
a) $\a1^2r_1l_1=1$,\\
b) $\car(-\un{k})=-r_1^{-k}l_1^{-k}\car(\un{k})$ for any $k\in\Z$,\\
c) for any $k\in\Z$ and $\fm\in G$ we have
\begin{equation}\label{eq:add}
\a1^k\car(\fm+\un{k})+\a1^{-k}\car(\fm-\un{k})=
\big((\a1r_1)^k+(\a1r_1)^{-k}\big)\car(\fm).
\end{equation}
\end{lem}
\begin{proof}
For a), multiply the two equations in \eqref{eq:thmHopfDe23} and apply
the multiplication map to both sides to obtain
\[\si(l_+l_-r_+r_-)=\a1^2l_+l_-r_+r_-.\]
Evaluate both sides at $\un{1}$ to get
\[1=lr(\un{0})=lr(\si^{-1}(\un{1}))=\si(lr)(\un{1})=\a1^2lr(\un{1})=
\a1^2l_1r_1.\]
Next \eqref{eq:hconds} gives for any $k\in\Z$,
\[0=\car(\un{k}-\un{k})=\car(\un{k})r_1^{-k}+l_1^k\car(-\un{k}),\]
hence b) follows.
Finally, using \eqref{eq:hconds} again, we have
\begin{align*}
\a1^k\car(\fm+\un{k})+\a1^{-k}\car(\fm-\un{k})&=
\a1^k\car(\fm)r_1^k+\a1^k l(\fm)\car(\un{k})+\\
&\quad +\a1^{-k}\car(\fm)r_1^{-k}+\a1^{-k}l(\fm)\car(-\un{k})=\\
&=\car(\fm)\big((\a1 r_1)^k+(\a1 r_1)^{-k}\big)+\\
&\quad +l(\fm)\car(\un{k})(\a1^k-\a1^{-k}r_1^{-k}l_1^{-k}).
\end{align*}
In the last equality we used part b). Now the second
term in the last expression vanishes due to part a). Thus c) follows.
\end{proof}

In what follows, we will treat the two cases when $\car(\un{1})=0$
and $\car(\un{1})\neq 0$ separately. The algebras satisfying the former
condition have a representation theory which reminds of that of the
enveloping algebra $U(\mathfrak{h}_3)$ of the three-dimensional
Heisenberg Lie algebra, while the latter case includes $U(\sl2)$ and
other algebras with similar structure of representations.

%
%
\subsubsection{The case $\car(\un{1})=0$.}
\begin{prop} \label{prop:h(1)=0}
If $\car(\un{1})=0$, then
\begin{equation}  \label{eq:h(1)=0}
\a1^2=r_1^2=1,\; \si(\car)=r_1\car,\; \si(r)=r_1 r,\;
\text{and }\si(l)= r_1 l.  
\end{equation}
In particular, $\langle X_+, X_-, \car\rangle $ is a subalgebra of $A$
with relations
\begin{align*}
  [X_+,X_-]&=h,  & [h,X_\pm]&=0,   &&&&\text{if $\xi=1$, $r_1=1$,}\\
  [X_+,X_-]&=h,  & \{h,X_\pm\}&=0, &&&&     \text{if $\xi=1$, $r_1=-1$,}\\
  \{X_+,X_-\}&=h,& [h,X_\pm]&=0,   &&&&       \text{if $\xi=-1$, $r_1=1$,}\\
  \{X_+,X_-\}&=h,& \{h,X_\pm\}&=0, &&&&    \text{if $\xi=-1$, $r_1=-1$,}
\end{align*}
respectively, where $\{\cdot,\cdot\}$ denotes anti-commutator.
\end{prop}
\begin{proof} Suppose $\car(\un{1})=0$. Then, by Lemma \ref{lem:dimabc}b),
$\car(-\un{1})=0$. This means that $\car\in\un{-1}=
\si^{-1}(\un{0})=\si^{-1}(\ker\ep)$. Thus $\ep(\si(\car))=0$.
Using \eqref{eq:HopfAxiomCounit}, \eqref{eq:thmHopfsiDe} and
\eqref{eq:thmHopfDe1} we deduce
\begin{align*}
\si(\car)&=(\ep\otimes 1)(\De(\si(\car)))=(\ep\otimes 1)(\si\otimes 1)(\De(\car))=\\
&=\ep(\si(\car))\otimes r+\ep(\si(l))\otimes \car=\ep(\si(l))\car.
\end{align*}
Analogously one proves $\si(\car)=\ep(\si(r))\car$.
Hence $\ep(\si(r))=\ep(\si(l))$.
But
\[\ep(\si(r))=\si(r) (\ker \ep)=\si(r)(\un{0})=r(-\un{1})=r_1^{-1}\]
and similarly for $l$. So $r_1=l_1$. From Lemma \ref{lem:dimabc}a)
we obtain $(\a1 r_1)^2=1$.
Now
\[S(\si(\car))=S(r_1^{-1}\car )=-r_1^{-1}\car,\quad\text{and}\quad
\si^{-1}(S(\car))=\si^{-1}(-h)=-r_1\car,\]
so \eqref{eq:TH_S1} implies that $r_1^2=1$. A similar calculation as above
shows that $\si(r)=r_1^{-1} r=r_1 r$ and $\si(l)=l_1^{-1} l= r_1 l$.
\end{proof}

We leave it to the reader to prove the following statement.
\begin{prop}
All finite-dimensional simple modules over an algebra $A(R,\si,\car,\xi)$
satisfying \eqref{eq:h(1)=0} and one of the commutation relations above
are either one- or two-dimensional.
\end{prop}
\begin{rem}
The algebra $U(\mathfrak{h}_3)$ is an ambiskew polynomial ring,
as shown in Section \ref{sec:U(h_3)}. For this algebra we have
$\car(\un{1})=0$ and $\xi=r_1=1$.
\end{rem}

%
%
\subsubsection{The case $\car(\un{1})\neq 0$.}
In this section, we consider the more complicated case when
$\car(\un{1})\neq 0$. We prove Theorem \ref{thm:mainresult_ny} which describes
the dimensions of $L(\fm)$ in terms of $\fm$.
The following two subsets of $G$ will play a vital role:
\begin{align}\label{eq:defG0}
G_0&=\{\fm\in G\; |\; \car(\fm)=0\},\\ \label{eq:defG1/2}
G_{1/2}&=\{\fm\in G\; |\; \car(\fm-\un{1})+\a1\car(\fm)=0\}.
\end{align}
The reason for this notation is that when $A=U(\sl2)$ as in 
Section \ref{sec:Usl2} then we have $G_0=\{(H-0)\}$ and
$G_{1/2}=\{(H-\frac{1}{2})\}$.
From \eqref{eq:hconds} it is immediate that $G_0$ is a subgroup
of $G$.
By Proposition \ref{prop:dimLm} we have
\begin{equation}
  \label{eq:G0dimLm1}
  G_0=\{\fm\in G\; |\; \dim L(\fm)=1\}.
\end{equation}
The following analogous result holds for $G_{1/2}$.
\begin{prop} \label{prop:dimLm2}
\begin{equation}  \label{eq:G1/2dimLm}
G_{1/2}=\{\fm\in G\; |\; \dim L(\fm)=2\}.
\end{equation}
\end{prop}
\begin{proof} If $\fm\in G_{1/2}$, then
by Proposition \ref{prop:dimLm}, $\dim L(\fm)\le 2$. But if
$\dim L(\fm)=1$, then $\car(\fm)=0$ so using $\fm\in G_{1/2}$ we get
$\car(\fm-\un{1})=0$ also. Since $G_0$ is a group we
deduce that $\un{1}\in G_0$, i.e. $\car(\un{1})=0$ which is
a contradiction. So $\dim L(\fm)=2$.
The converse inclusion is immediate from Proposition \ref{prop:dimLm}.
\end{proof}

Set
\begin{equation}
  \label{eq:Ndef}
  N=\begin{cases}\text{order of $\xi r_1$}&\text{if $(\xi r_1)^2\neq 1$ and
$\xi r_1$ is a root of unity,}\\
\infty&\text{otherwise.}\end{cases}
\end{equation}
We also set
\[N'=\begin{cases}N,&\text{if }N\text{ is odd},\\
N/2,&\text{if }N\text{ is even},\\
\infty,&\text{if }N=\infty.\end{cases}\]

The next statement describes the intersection of $G_0$
and $G_{1/2}$ with $\un{\Z}$.
\begin{prop}\label{prop:G_0nZ}
We have
\begin{equation}\label{eq:G_0nZ}
G_0\cap \un{\Z}=
\begin{cases}
\{\un{0}\},&\text{if }N=\infty,\\
\un{N'\Z},&\text{otherwise},
\end{cases}
\end{equation}
and
\begin{equation}\label{eq:G12nZ}
G_{1/2}\cap \un{\Z_{>0}}=
\begin{cases}
\emptyset,&\text{if }N=\infty,\\
\{\un{n}\in\un{\Z_{>0}} : N\big| 2n-1\},&\text{otherwise}.
\end{cases}
\end{equation}
\end{prop}
\begin{rem} The set $G_{1/2}\cap\un{\Z_{\le 0}}$ can be understood
using \eqref{eq:G12nZ} and Lemma \ref{lem:dim2special}a).
\end{rem}
\begin{proof}
We first prove \eqref{eq:G_0nZ}. Let $n\in\Z$. The right hand side
of \eqref{eq:G_0nZ} is invariant under $n\mapsto -n$. By Lemma
\ref{lem:dimabc}b) so is the left hand side. Moreover since
$\car(\un{0})=0$, the ideal $\un{0}$ belongs to both sides of the
equality. Thus we can assume $n>0$.

Using \eqref{eq:hconds} and that $r$ and $l$, viewed as functions $G\to\K$,
are multiplicative homomorphisms it follows by induction that
\[\car(\un{n})=\car(\un{1})\sum_{i=0}^{n-1}r_1^il_1^{n-1-i}.\]
By Lemma \ref{lem:dimabc}a),
$r_1/l_1=(\a1 r_1)^2/(\a1^2 r_1 l_1)=(\a1 r_1)^2$, so we can rewrite this as
\begin{equation}\label{eq:hnformula}
\car(\un{n})=\car(\un{1})l_1^{n-1}\sum_{i=0}^{n-1} (\a1 r_1)^{2i}.
\end{equation}
If $N=\infty$ and $(\xi r_1)^2\neq 1$ then by \eqref{eq:hnformula} we have
$\un{n}\in G_0\cap\un{\Z}$ iff $(\xi r_1)^{2n}=1$, which is false.
If $(\xi r_1)^2=1$, then \eqref{eq:hnformula} implies
that $\un{n}\notin G_0\cap\un{\Z}$.
If $N<\infty$, then $(\xi r_1)^2\neq 1$ so by
\eqref{eq:hnformula}, $\car(\un{n})=0$ iff $(\xi r_1)^{2n}=1$
i.e. iff $N\big| 2n$. This is equivalent to $N'\big| n$.

Next we prove \eqref{eq:G12nZ}. Suppose $n\in\Z_{>0}$. By
definition, $n\in G_{1/2}$ iff
\[\car(\un{n}-\un{1})+\xi\car(\un{n})=0.\] Using
\eqref{eq:hnformula} on both terms and dividing by
$\car(\un{1})\xi l_1^{n-1}$, this is equivalent to
\[\xi^{-1}l_1^{-1}\sum_{k=0}^{n-2} (\xi r_1)^{2k}+\sum_{k=0}^{n-1}(\xi r_1)^{2k}=0. \]
But $\xi^{-1}l_1^{-1}=\xi r_1$ by Lemma \ref{lem:dimabc}a) so this
can be rewritten as
\begin{equation}\label{eq:dim12hej}
\sum_{k=0}^{2n-2} (\xi r_1)^k=0.
\end{equation}
Thus $(\xi r_1)^2\neq 1$ and multiplying by $\xi r_1-1$ we get
$(\xi r_1)^{2n-1}=1$. Therefore $N<\infty$ and $N\big| 2n-1$.
Conversely, if $N<\infty$ and $N\big| 2n-1$ then $(\xi r_1)^2\neq
1$ and $(\xi r_1)^{2n-1}=1$ which implies \eqref{eq:dim12hej}. This
proves \eqref{eq:G12nZ}.
\end{proof}


\begin{prop} \label{prop:G12coset}
Suppose $\car(\un{1})\neq 0$ and $G_{1/2}\neq\emptyset$. Then\\
a) $\xi r_1\neq -1$, and\\
b) $G_{1/2}$ is a left and right coset of $G_0$ in $G$.
\end{prop}
\begin{proof}
Let $\fm_{1/2}\in G_{1/2}$.
To prove a), suppose that $\xi r_1=-1$. Then
\begin{align*}
0&=\car(\fm_{1/2}-\un{1})+\xi\car(\fm_{1/2})=\\&=
\car(\fm_{1/2})r(-\un{1})+l(\fm_{1/2})\car(-\un{1})+\xi\car(\fm_{1/2})=\\
&=\car(\fm_{1/2})(r_1^{-1}+\xi)+l(\fm_{1/2})\car(-\un{1})=\\
&=-l(\fm_{1/2})r_1^{-1}l_1^{-1}\car(\un{1}),
\end{align*}
where we used Lemma \ref{lem:dimabc}b) in the last equality.
Since $l$ is
invertible we deduce that $\car(\un{1})=0$ which is a contradiction.

To prove part b), we will show that
\[G_{1/2}=G_0+\fm_{1/2}.\]
One proves $G_{1/2}=\fm_{1/2}+G_0$ in an analogous way.
Let $\fm\in G_0$ be arbitrary.
Then using \eqref{eq:hconds} twice,
\[\car(\fm+\fm_{1/2}-\un{1})+\xi\car(\fm+\fm_{1/2})=
l(\fm)\big(\car(\fm_{1/2}-\un{1})+\xi\car(\fm_{1/2})\big)=0.\]
Since $l$ is invertible we get $\fm+\fm_{1/2}\in G_{1/2}$.

Conversely, suppose $\fm\in G_{1/2}$. Then
\begin{align*}
\car(\fm-\un{1})+\xi\car(\fm)&=0,\\
\car(\fm_{1/2}-\un{1})+\xi\car(\fm_{1/2})&=0.
\end{align*}
Multiply the first equation by $r(-\fm_{1/2})$ and the second
by $-r(-\fm_{1/2})l(-\fm_{1/2})l(\fm)$ and add them together.
Then we get
\begin{multline*}
\big((\car(\fm)r_1^{-1}+l(\fm)\car(-\un{1})\big)r(-\fm_{1/2})-\\
r(-\fm_{1/2})l(-\fm_{1/2})l(\fm)\big(\car(\fm_{1/2}r_1^{-1}+
l(\fm_{1/2})\car(-\un{1})\big)+\xi\car(\fm-\fm_{1/2})=0,
\end{multline*}
or equivalently,
\[\car(\fm)r_1^{-1}r(-\fm_{1/2})-r(-\fm_{1/2})l(-\fm_{1/2})l(\fm)
\car(\fm_{1/2})r_1^{-1}+\xi\car(\fm-\fm_{1/2})=0.\]
Using \eqref{eq:hconds} this can be written
\[r_1^{-1}(1+\xi r_1)\car(\fm-\fm_{1/2})=0.\]
Since $\xi r_1\neq -1$ by part a), we conclude that
$\car(\fm-\fm_{1/2})=0$. This shows that $\fm\in G_0+\fm_{1/2}$.
\end{proof}

The following lemma will be useful.
\begin{lem}\label{lem:atmost}
Let $j\in\Z$. If $\fm_0\in G_0$, then
\begin{equation}
  \label{eq:G0iff}
\fm_0+\un{j}\in G_0 \Longleftrightarrow \un{j}\in G_0,
\end{equation}
and if $\car(\un{1})\neq 0$ and $\fm_{1/2}\in G_{1/2}$, then
\begin{equation}
  \label{eq:G12iff}
\fm_{1/2}+\un{j}\in G_{1/2} \Longleftrightarrow \un{j}\in G_0.
\end{equation}
\end{lem}
\begin{proof}
\eqref{eq:G0iff} is immediate since $G_0$ is a subgroup of $G$.
If $\un{j}\in G_0$, then $\fm_{1/2}+\un{j}\in G_{1/2}$ by
Proposition \ref{prop:G12coset}.
Conversely, if $\fm_{1/2}+\un{j}\in G_{1/2}$
then by Proposition \ref{prop:G12coset},
$G_0\ni \fm_{1/2}+\un{j}-\fm_{1/2}=\un{j}$.
%
\end{proof}

%
%

The next statements will
be needed in Section \ref{sec:semisimplicity}.

\begin{lem}\label{lem:dim2special}
Suppose $\car(\un{1})\neq 0$ and let $\fm,\fn \in G_{1/2}$.
Then 
\begin{itemize}
\item[a)] $\un{1}-\fm\in G_{1/2}$, and
\item[b)] $\fm+\fn-\un{1} \in G_0$.
\end{itemize}
\end{lem}
\begin{proof}
Part a) follows from the calculation
\begin{align*}
\car(\un{1}-\fm-\un{1})+\xi\car(\un{1}-\fm)&=
-l(-\fm)r(-\fm)\car(\fm)-\xi l(\un{1}-\fm)(r(\un{1}-\fm)\car(\fm-\un{1})=\\
&=-l(-\fm)r(-\fm)\big(\car(\fm)+\xi r_1l_1\car(\fm-\un{1})\big)=\\
&=-l(-\fm)r(-\fm)\xi^{-1}\big(\xi\car(\fm)+ \car(\fm-\un{1})\big)=0.
\end{align*}
For part b), use that $\dim L(\un{1}-\fn)=2$ by part a), and thus
$\fm+\fn-\un{1}=\fm-(\un{1}-\fn)\in G_0$ by Proposition \ref{prop:G12coset}b).
\end{proof}

The formulas provided by the following technical lemma are the
key to proving our main theorem.
\begin{lem} Let $\fm\in G$ and $j\in\Z_{\ge 0}$. If $n=2j+1$ then
  \begin{equation}
    \label{eq:oddn}
    \sum_{k=0}^{n-1}\xi^{n-1-k}\car(\fm-\un{k})=
r_1^{-j}\car(\fm-\un{j})\sum_{k=0}^{n-1} (\xi r_1)^k
  \end{equation}
and if $n=2j+2$ then
\begin{equation}
  \label{eq:evenn}
   \sum_{k=0}^{n-1}\xi^{n-1-k}\car(\fm-\un{k})
  =r_1^{-j}\big(\car(\fm-\un{j}-\un{1})+\xi\car(\fm-\un{j})\big)
\sum_{k=0}^{n/2-1}(\xi r_1)^{2k}.
\end{equation}
\end{lem}
\begin{proof}
If $n=2j+1$,
we make the change of index $k\mapsto j-k$, then
factor out $\a1^j$ and apply formula \eqref{eq:add}:
\[\sum_{k=0}^{2j}\xi^{2j-k}\car(\fm-\un{k})=
\sum_{k=-j}^j\xi^{j+k}\car(\fm-\un{j}+\un{k})=
\xi^j \car(\fm-\un{j})\sum_{k=-j}^j (\xi r_1)^k.\]
Factoring out $(\xi r_1)^{-j}$ and changing index from $k$ to $k-j$
yields \eqref{eq:oddn}.

For the $n=2j+2$ case we first split the sum in the left hand side
of \eqref{eq:evenn} into two sums corresponding to odd and even $k$:
\[\sum_{k=0}^j\xi^{2j-2k}\car(\fm-\un{2k}-\un{1})+
\sum_{k=0}^j\xi^{2j+1-2k}\car(\fm-\un{2k})\]
Then we make the change of summation index $k\mapsto -k+j/2$ in both sums
\[
\xi^j\sum_{k=-j/2}^{j/2}\xi^{2k}\car(\fm-\un{j}-\un{1}+\un{2k})+
\xi^{j+1}\sum_{k=-j/2}^{j/2}\xi^{2k}\car(\fm-\un{j}+\un{2k})
\]
and use \eqref{eq:add} on each of them to get
\[\big(\car(\fm-\un{j}-\un{1})+\xi\car(\fm-\un{j})\big)
\xi^j\sum_{k=-j/2}^{j/2}(\xi r_1)^{2k}.\]
If we factor out $(\xi r_1)^{-j}$ and change summation index
from $k$ to $k-j/2$ we obtain \eqref{eq:evenn}.
\end{proof}

%
%

We now come to the main results in this section.
\begin{mainlem} \label{lem:mainresult}
Assume that $\car(\un{1})\neq 0$ and let $\fm\in G$. Then
\begin{itemize}
\item[a)] $\dim L(\fm)\le N$,
\item[b)] if $\dim L(\fm)=n<N$ then $\fm\in G_{\frac{i-1}{2}}+\un{j}$ where
$n=2j+i$, $i\in\{1,2\}$, $j\in\Z_{\ge 0}$, and
\item[c)] if $i\in\{1,2\}$, $j\in\Z_{\ge 0}$, $2j+i\le N$ and $\fm\in G_{\frac{i-1}{2}}$ then
\begin{equation}  \label{eq:mainthm1}
\dim L(\fm+\un{j})=2j+i.
\end{equation}
\item[d)] If $N'<\infty$ then $\dim L(\fm+\un{N'j})=\dim L(\fm)$ for any
$j\in\Z$.
\end{itemize}
\end{mainlem}
\begin{proof}
Part a) is trivial when $N=\infty$. If $N$ is finite and odd,
Proposition \ref{prop:dimLm} and \eqref{eq:oddn} imply that $\dim L(\fm)\le N$.
If $N$ is finite and even, then $(\xi r_1)^N=1$ and $(\xi r_1)^2\neq 1$
so $\sum_{k=0}^{N/2-1}(\xi r_1)^{2k}=0$. Hence Proposition \ref{prop:dimLm}
and \eqref{eq:evenn} implies $\dim L(\fm)\le N$ in this case as well.

Next we turn to part b). Suppose first that $\dim L(\fm)=n=2j+1<N$.
Then by Proposition \ref{prop:dimLm} and \eqref{eq:oddn} the right
hand side of \eqref{eq:oddn} is zero. The definition of $N$ implies
that $\car(\fm-\un{j})=0$, i.e. $\fm\in G_0+\un{j}$.
If instead $\dim L(\fm)=2j+2<N$, Proposition \ref{prop:dimLm}
and \eqref{eq:evenn}
similarly implies that $\fm\in G_{1/2}+\un{j}$.

To prove \eqref{eq:mainthm1},
we proceed by induction on $j$. For $j=0$ it follows
from \eqref{eq:G0dimLm1} and \eqref{eq:G1/2dimLm}. Suppose it
holds for $j=0,1,\ldots,k-1$, where $k>0$ and $2k+i\le N$.
We first show that $\dim L(\fm+\un{k})\le 2k+i$.
If $i=1$ then by \eqref{eq:oddn},
\[\sum_{l=0}^{2k}\xi^{2k-l}\car(\fm+\un{k}-\un{l})=
r_1^{-k}\car(\fm)\sum_{l=0}^{2k} (\xi r_1)^l=0\]
since $\fm\in G_0$. Similarly, if $i=2$, then \eqref{eq:evenn} gives
\[\sum_{l=0}^{2k+1}\xi^{2k+1-l}\car(\fm+\un{k}-\un{l})
  =r_1^{-k}\big(\car(\fm-\un{1})+\xi\car(\fm)\big)
\sum_{l=0}^k(\xi r_1)^{2l}=0\]
since $\fm\in G_{1/2}$ in this case.
Thus $\dim L(\fm+\un{j})\le 2j+i$ by Proposition \ref{prop:dimLm}.
Write $\dim L(\fm+\un{k})=2k'+i'$ where $k'\ge 0$, $i'\in\{1,2\}$ and
assume that $2k'+i'<2k+i$. By part b)
we have $\fm+\un{k}\in G_{\frac{i'-1}{2}}+\un{k'}$ which implies
that $\dim L(\fm+\un{k}-\un{k'})=i'$ by \eqref{eq:G0dimLm1} and
\eqref{eq:G1/2dimLm}. This contradicts the induction hypothesis
unless $k'=0$. Assuming $k'=0$ we get $\fm+\un{k}\in
G_{\frac{i'-1}{2}}$. If $i=i'$
then from Lemma \ref{lem:atmost} follows that $\un{k}\in G_0$.
Since $0<k<\frac{2k+i}{2}\le N/2\le N'$ this contradicts
\ref{eq:G_0nZ}.
We now show that $i\neq i'$ is also impossible.
If $i=1$ and $i'=2$, then $\fm\in G_0$ and $\fm+\un{k}\in G_{1/2}$
so by Proposition \ref{prop:G12coset}b), $\un{k}\in G_{1/2}\cap
\un{\Z_{>0}}$. By \eqref{eq:G12nZ} we get $N\big| 2k-1$ which is
absurd because $0<2k-1<2k+1\le N$. If $i=2$ and $i'=1$ then
$\fm\in G_{1/2}$ and $\fm+\un{k}\in G_0$. By Proposition \ref{prop:G12coset}b)
we have $\un{-k}=\fm-(\fm+\un{k})\in
G_{1/2}$. By Lemma \ref{lem:dim2special}a), $\un{1+k}\in G_{1/2}$
so \eqref{eq:G12nZ} implies that $N\big| 2(1+k)-1=2k+1$. This
is impossible since $0<2k+1<2k+2\le N$.
We have proved that the assumption $2k'+i'<2k+i$ is false and
hence that $\dim L(\fm+\un{k})=2k+i$, which proves the induction
step.

Finally, part d) follows from Corollary \ref{cor:period} and
Proposition \ref{prop:G_0nZ}.
\end{proof}

\begin{thm}\label{thm:mainresult_ny}
Let $\fm\in G$. 
\begin{itemize}
\item If $N=\infty$, then
\begin{gather}
\label{eq:MRny1}
\dim L(\fm)<\infty \Longleftrightarrow \fm\in (G_0+\un{\Z_{\ge 0}})
\cup (G_{1/2}+\un{\Z_{\ge 0}})\\
\intertext{and}
\label{eq:MRny2}
\dim L(\fm_0+\un{j})=2j+1,\quad\text{for $\fm_0\in G_0$ and
$j\in\Z_{\ge 0}$},\\
\label{eq:MRny3}
\dim L(\fm_{1/2}+\un{j})=2j+2,\quad\text{for $\fm_{1/2}\in G_{1/2}$ and
$j\in\Z_{\ge 0}$}.
\end{gather}
\item If $N<\infty$ and $N$ is even, then
\begin{gather}
\label{eq:MRny4}
\dim L(\fm)<\infty \Longleftrightarrow \fm\in (G_0+\un{\Z})
\cup (G_{1/2}+\un{\Z})\\
\intertext{and}
\label{eq:MRny5}
\dim L(\fm+\un{(N/2)j})=\dim L(\fm),\quad\text{for any 
$\fm\in G$ and $j\in\Z$},\\
\intertext{and for $\fm_0\in G_0$ and $\fm_{1/2}\in G_{1/2}$ we have}
\label{eq:MRny6}
\dim L(\fm_0+\un{j})=2j+1, \quad\text{if $0\le j<N/2$},\\
\label{eq:MRny7}
\dim L(\fm_{1/2}+\un{j})=2j+2,
\quad\text{if $0\le j< N/2$}.
\end{gather}
\item If $N<\infty$ and $N$ is odd, then
\begin{gather}
\label{eq:MRny8}
\dim L(\fm)<\infty \Longleftrightarrow \fm\in G_0+\un{\Z}=G_{1/2}+\un{\Z}\\
\intertext{and}
\label{eq:MRny9}
\dim L(\fm+\un{Nj})=\dim L(\fm),\quad
\text{for any $\fm\in G$ and $j\in\Z$},\\
\intertext{and for $\fm_0\in G_0$ and $\fm_{1/2}\in G_{1/2}$ we have}
\label{eq:MRny10}
\dim L(\fm_0+\un{j})=\begin{cases}2j+1,&\text{if $0\le j< \frac{N+1}{2}$},\\
2j+1-N,&\text{if $\frac{N+1}{2}\le j <N$,}
\end{cases}\\
\label{eq:MRny11}
\dim L(\fm_{1/2}+\un{j})=
\begin{cases}2j+2,&\text{if $0\le j< \frac{N-1}{2}$},\\
2j+2-N,&\text{if $\frac{N-1}{2}\le j< N$.}
\end{cases}
\end{gather}
\end{itemize}
\end{thm}
\begin{proof}
When $N=\infty$, relations \eqref{eq:MRny1}-\eqref{eq:MRny3} are immediate
from Lemma \ref{lem:mainresult}b) and c).

Suppose $N$ is finite and even.
The $\Rightarrow$ implication in \eqref{eq:MRny4} holds
by Lemma \ref{lem:mainresult}b). And
\eqref{eq:MRny5} follows from \eqref{eq:G_0nZ} and
Corollary \ref{cor:period}. Assume that $\fm\in (G_0+\un{\Z})\cup
(G_{1/2}+\un{\Z})$. Using \eqref{eq:MRny5} we can assume that
$\fm = \fm'+\un{j}$ where $\fm'\in G_0\cup G_{1/2}$ and
$0\le j<N/2$. Then, if $i\in\{1,2\}$ we have $2j+i\le N$ and
Lemma \ref{lem:mainresult}c) implies 
\eqref{eq:MRny6}-\eqref{eq:MRny7} and therefore $\dim L(\fm)<\infty$
so \eqref{eq:MRny4} is also proved.

Assume that $N$ is finite and odd.
By \eqref{eq:G12nZ} we have $\un{(N+1)/2}\in G_{1/2}$.
Therefore $G_0+\un{\Z}=G_0+\un{(N+1)/2}+\un{\Z}=G_{1/2}+\un{\Z}$
since $G_{1/2}$ is a right coset of $G_0$ in $G$ by
Proposition \ref{prop:G12coset}. As before,
Lemma \ref{lem:mainresult}b) implies the $\Rightarrow$ case in
\eqref{eq:MRny8} and \eqref{eq:MRny9} holds by virtue
of \eqref{eq:G_0nZ} and Corollary \ref{cor:period}.
If $\fm\in G_0+\un{\Z}$ we can assume by \eqref{eq:MRny9} that
$\fm\in G_0+\un{j}$ where $0\le j< N$. If $j<\frac{N+1}{2}$, then
$2j+1<N+2$ so since $N$ is odd we have $2j+1\le N$. By
Lemma \ref{lem:mainresult}c) we deduce that $\dim L(\fm)=2j+1$.
If instead $j\ge\frac{N+1}{2}$, then
$\fm=\un{(N+1)/2}+\fm-\un{(N+1)/2}\in G_{1/2}+\un{k}$
where $k=j-\frac{N+1}{2}$ so $0\le k<\frac{N-1}{2}$.
Thus $2k+2\le N$ so Lemma \ref{lem:mainresult}c)
implies that $\dim L(\fm)=2k+2=2j+1-N$.
This proves \eqref{eq:MRny10} and the $\Leftarrow$ implication
in \eqref{eq:MRny8}.
Finally \eqref{eq:MRny11} is equivalent to \eqref{eq:MRny10}
in the following sense. Let $0\le j<N$ and $\fm_{1/2}\in G_{1/2}$.
Then \[\dim L(\fm_{1/2}+\un{j})=\dim L(\fm_0+\un{j'}),\]
where
$j'=j+(N+1)/2$ and $\fm_0=\fm_{1/2}-\un{(N+1)/2}$. Now $\fm_0\in G_0$ since
$G_{1/2}$ is a coset of $G_0$ in $G$. If $0\le j<\frac{N-1}{2}$,
then $\frac{N+1}{2}\le j'<N$ so by \eqref{eq:MRny10} we have
\[\dim L(\fm_{1/2}+\un{j})=
\dim L(\fm_0+\un{j'})=2j'+1-N=2j+2.\]
And if $\frac{N-1}{2}\le j< N$,
then $0\le j'-N<\frac{N+1}{2}$ and hence
\[\dim L(\fm_{1/2}+\un{j})=
\dim L(\fm_0+\un{j'-N})=2(j'-N)+1=2j+1-N.\]
The proof is finished.
\end{proof}

\begin{cor}\label{cor:infdim}
If $N=\infty$ and $\fm\in G_0\cup G_{1/2}$, then
$L(\fm+\un{j})$ is infinite-dimensional for any $j\in\Z_{<0}$.
\end{cor}
\begin{proof}
If the dimension of $L(\fm+\un{j})$ were finite and odd (even),
then $\dim L(\fm+\un{j-k})=1\; (2)$ for some
$k\ge 0$ by Lemma \ref{lem:mainresult}b). By Lemma \ref{lem:mainresult}c),
$L(\fm)$ has then dimension $2(j-k)+1$ ($2(j-k)+2$) and thus $j=k$
which is absurd.
\end{proof}

\begin{cor}\label{cor:vermastruct}
Suppose $N=\infty$ and
let $\fm\in G_f$. Then $L(\fm)$ is the unique finite-dimensional
quotient of $M(\fm)$.
\end{cor}
\begin{proof} It is enough to prove that the unique maximal proper
submodule $N(\fm)$ of $M(\fm)$ is simple.
By Theorem \ref{thm:mainresult_ny} we can write $\fm=\fn+\un{j}$
where $\fn\in G_0\cup G_{1/2}$ and $j\in\Z_{\ge 0}$.
From the proof of Proposition \ref{prop:dimLm} we have
\[\supp(L(\fm))=\{\fn+\un{j},\fn+\un{j}-\un{1},\ldots,\fn-\un{j}\}.\]
Thus $N(\fm)$ is a highest weight module of highest weight
$\fn-\un{j}-\un{1}$. So $N(\fm)$ is a quotient of
$M(\fn-\un{j}-\un{1})$. But $M(\fn-\un{j}-\un{1})$ is simple,
otherwise it would have a finite-dimensional simple quotient,
i.e. $L(\fn-\un{j}-\un{1})$ would be finite-dimensional,
contradicting Corollary \ref{cor:infdim}. Thus $N(\fm)$ is
also simple.
\end{proof}

\begin{rem}
We finish this section by remarking that there exist algebras
in the class studied in this paper which do not have even-dimensional
simple modules as for example the algebra $B_q$ from Section \ref{sec:downup}.
Indeed, in this case we have
$\xi r_1=-1$ and so $N=\infty$ by definition.
By Proposition \ref{prop:G12coset}, $G_{1/2}=\emptyset$
so by Theorem \ref{thm:mainresult_ny}, there can exist
no even-dimensional simple modules.
\end{rem}

%
\section{Tensor products and a Clebsch-Gordan formula}
%
\label{sec:CG}
As we have seen in Section \ref{sec:prel}
the existence of a Hopf structure on an algebra allows
one to define tensor product of its representations by
\eqref{eq:tensordef}.
The aim of this section is to prove
a formula which decomposes the tensor product of
two simple $A$-modules into a direct sum of simple modules.
It generalizes the classical Clebsch-Gordan formula
for modules over $U(\sl2)$. We will assume that $A=A(R,\si,\car,\xi)$
is an ambiskew polynomial ring and that it carries a Hopf
structure of the type considered in Section \ref{sec:hopf}.
We will also assume \eqref{eq:torfree} and that $N=\infty$.

\begin{lem}\label{lem:tensorweight}
Let $V$ and $W$ be two $A$-modules. Then
\begin{equation}\label{eq:proptensor1}
V_\fm\otimes W_\fn \subseteq (V\otimes W)_{\fm+\fn}
\end{equation}
for any $\fm,\fn\in G$.
Hence if $V$ and $W$ are weight modules, then so is
$V\otimes W$ and
\[\supp(V\otimes W)=\{\fm+\fn \, |\, \fm\in\supp(V),\fn\in\supp(W)\}.\]
\end{lem}
\begin{proof} Let $v\in V_\fm, w\in W_\fn$. Then for any $r\in R$,
\begin{multline*}
r(v\otimes w)=\sum_{(r)}r'v\otimes r''w =
\sum_{(r)}r'(\fm)v\otimes r''(\fn)w =\\=
\sum_{(r)}r'(\fm)r''(\fn) v\otimes w = r(\fm+\fn) v\otimes w
\end{multline*}
by \eqref{eq:grplem2}, proving \eqref{eq:proptensor1}. Thus if
$V,W$ are weight modules,
\[V\otimes W=(\oplus_\fm V_\fm)\otimes(\oplus_\fn W_\fn)=
\oplus_{\fm,\fn}V_\fm\otimes W_\fn = \oplus_\fm\big(\oplus_{\fm_1+\fm_2=\fm}
V_{\fm_1}\otimes W_{\fm_2}\big).\]
\end{proof}

\begin{thm} Let $\fm,\fn\in G_f$.
We have the following isomorphism
\begin{equation}\label{eq:clebschgordan}
L(\fm)\otimes L(\fn)\simeq L(\fm+\fn)\oplus
L(\fm+\fn-\un{1})\oplus\ldots\oplus L(\fm+\fn-\un{s}+\un{1})
\end{equation}
where $s=\min\{\dim L(\fm),\dim L(\fn)\}$.
\end{thm}
\begin{proof}
Let $e^\fm, e^\fn$ denote highest weight vectors in $L(\fm)$,
$L(\fn)$ respectively
and set $e^\fm_j:=(X_-)^je^\fm$ for $j\in\Z_{\ge 0}$
and similarly for $\fn$. Set $V=L(\fm)\otimes L(\fn)$.
By Lemma \ref{lem:tensorweight} we have
\[V_{\fm+\fn-\un{k}}=\oplus_{i+j=k} \K e^\fm_i\otimes e^\fn_j\]
for $k\in\Z_{\ge 0}$.
Fix $0\le k\le s-1$.
We will prove that
\begin{equation}\label{eq:TCdimeq}
\dim \ker X_+|_{V_{\fm+\fn-\un{k}}}=1.
\end{equation}
From the calculations in the proof of Proposition \ref{prop:dimLm}
follows that when $j>0$, $X_+e^\fm_j$ is a nonzero multiple of $e^\fm_{j-1}$.
Let $\nu_j^\fm$ denote this multiple.
Let
\[u=\sum_{i=0}^k\la_i e^\fm_i\otimes e^\fn_{k-i}\]
be an arbitrary vector in $V_{\fm+\fn-\un{k}}$.
Then
\begin{align*}
X_+u&=\sum_{i=0}^k \la_i(X_+e^\fm_i\otimes r_+e^\fn_{k-i}+
l_+e^\fm_i\otimes X_+e^\fn_{k-i})=\\
&=\sum_{i=0}^{k-1}\big[\la_{i+1}\nu^\fm_{i+1}r_+(\fn-\un{k}+\un{i}+\un{1})+
\la_il_+(\fm-\un{i})\nu^\fn_{k-i}\big]e^\fm_i\otimes e^\fn_{k-1-i}.
\end{align*}
Setting
\begin{align*}
c_i&=l_+(\fm-\un{i})\nu^\fn_{k-i},\\
c_i'&=\nu^\fm_ir_+(\fn-\un{k}+\un{i}),
\end{align*}
the condition for $u$ to be a highest weight vector can hence be written as

\begin{equation}\label{eq:TC_matrix}
\begin{bmatrix}
  c_0 & c_1'&       &   & \\
      & c_1 & c_2'  &    & \\
      &     & \ddots & \ddots &\\
 &&& c_{k-1} & c_k'
\end{bmatrix}
\begin{bmatrix}
  \la_0\\
  \la_1\\
  \vdots\\
  \la_k
\end{bmatrix}=0.
\end{equation}
Since $r_+$ and $l_+$ are grouplike, they are invertible and hence
$c_i\neq 0\neq c_{i+1}'$ for any $i=0,1,\ldots, k-1$. Therefore the
space of solutions to \eqref{eq:TC_matrix} is one-dimensional.
Thus \eqref{eq:TCdimeq} is proved.

From the definition of Verma modules, it follows that for
$k=0,1,\ldots, s-1$, there is a nonzero $A$-module morphism
\[M(\fm+\fn-\un{k})\to L(\fm)\otimes L(\fn)\]
which maps a highest weight vector in $M(\fm+\fn-\un{k})$ to
a highest weight vector in $L(\fm)\otimes L(\fn)$ of weight
$\fm+\fn-\un{k}$.
But $L(\fm)\otimes L(\fn)$ is finite-dimensional so this
morphism must factor through $L(\fm+\fn-\un{k})$ by Corollary
\ref{cor:vermastruct}. Taking direct sums of these morphisms
we obtain an $A$-module morphism
\[\ph:L(\fm+\fn)\oplus
L(\fm+\fn-\un{1})\oplus\ldots\oplus L(\fm+\fn-\un{s}+\un{1})
\to L(\fm)\otimes L(\fn).\]
We claim it is injective. Indeed,
the projection of the kernel of $\ph$ to any term
$L(\fm+\fn-\un{i})$ must be zero,
because it is a proper submodule of the simple module $L(\fm+\fn-\un{i})$.

To conclude
we now calculate the dimensions of both sides.
Write $\dim L(\fm)=2j_1+i_1$ and $\dim L(\fn)=2j_2+i_2$
where $j_1,j_2\in\Z_{\ge 0}$ and $i_1,i_2\in\{1,2\}$.
By Lemma~\ref{lem:mainresult}b),
$\dim L(\fm-\un{j_1})=i_1$ and $\dim L(\fn-\un{j_2})=i_2$.
First note that
\[\dim L(\fm-\un{j_1}+\fn-\un{j_2})=i_1+i_2-1.\]
When $i_1=i_2=1$, this is true because $G_0$ is a subgroup of $G$.
When one of $i_1,i_2$ is $1$ and the other $2$, it follows from
Proposition \ref{prop:G12coset}b). And if $i_1=i_2=2$, it follows from
Lemma \ref{lem:dim2special}b) and Theorem \ref{thm:mainresult_ny}.

From Theorem \ref{thm:mainresult_ny} also
follows that $\dim L(\fm+\un{k})=\dim L(\fm)+2k$ if $\dim L(\fm)<\infty$
and $k\in\Z_{\ge 0}$. Hence, recalling that
$s=\min\{\dim L(\fm),\dim L(\fn)\}$, we have
\begin{align*}
\sum_{k=0}^{s-1} \dim L(\fm+\fn-\un{k})&=\sum_{k=0}^{s-1}
\dim L(\fm-\un{j_1}+\fn-\un{j_2}+\un{j_1}+\un{j_2}-\un{k})=\\
&=\sum_{k=0}^{s-1}\big(i_1+i_2-1+2(j_1+j_2-k)\big)=\\
&=s(i_1+i_2-1+2j_1+2j_2)-s(s-1)=\\
&=s(\dim L(\fm)+\dim L(\fn)-s)=\\
&=\dim L(\fm) \dim L(\fn)=\dim \big(L(\fm)\otimes L(\fn)\big).
\end{align*}
This completes the proof of the theorem.
\end{proof}

Under some conditions it is possible to introduce a
$\ast$-structure on $A$. In this connection it would be interesting
to study Clebsch-Gordan coefficients and the relation with
special functions. This will be a subject for future investigation.

\section{Casimir operators and semisimplicity}
\label{sec:semisimplicity}
Arguing as in the proof of Lemma \ref{lem:Bfdweight}, it is easy to see
that any finite-dimensional
semisimple module over $A=A(R,\si,\car,\a1)$ is a weight module.
In this section we will prove the converse, that any finite-dimensional weight
module over $A$ is semisimple.
Note that in general not all finite-dimensional modules over our
algebra $A$ are semisimple. The corresponding example is constructed
in \cite{DQS} for the algebra from Section \ref{sec:dqs}.
A necessary and sufficient condition for all finite-dimensional
modules over an ambiskew polynomial ring to be semisimple was given
in \cite{Jor1}, Theorem 5.1.

In this section we assume that $A=A(R,\si,\car,\xi)$ is
an ambiskew polynomial ring with a Hopf structure of the type
introduced in Section \ref{sec:hopf} such that
\eqref{eq:torfree} holds. We also assume that $N=\infty$.

Let $V$ be a finite-dimensional weight module over $A$.
We will first treat the case when $\supp(V)\subseteq \fm+\un{\Z}$
where $\fm\in G_0$ is fixed.
Define a linear map
\[ C_V:V\to V\]
by requiring
\[ C_Vv = \si^j(t)v,\quad\text{for $v\in V_{\fm+\un{j}}$ and $j\in\Z$.}\]
Here $\si$ denotes the extended automorphism \eqref{eq:extaut}.
More explicitly we have (if $j\ge 0$)
\[C_Vv=\si^j(t)v=\Big(\xi^jt+\sum_{k=0}^{j-1}\xi^k\si^{j-1-k}(\car)\Big)v=
\xi^jtv+\sum_{k=0}^{j-1}\xi^k\car(\fm+\un{k+1})v\]
and similarly when $j<0$.
It is easy to check that $C_V$ is a morphism of $A$-modules.
Hence it is constant on each finite-dimensional simple module $V$
by Schur's Lemma. Moreover if $\ph:V\to W$ is a morphism of weight $A$-modules
with support in $\fm+\un{\Z}$, then $\ph C_V=C_W\ph $.
\begin{prop}\label{prop:separation}
Let $j_1,j_2\in\Z_{\ge 0}$. If
$C_{L(\fm+\un{j_1})}=C_{L(\fm+\un{j_2})}$, then $j_1=j_2$.
\end{prop}
\begin{proof}
By applying $C_{L(\fm+\un{j})}$ to the highest weight vector
of $L(\fm+\un{j})$, ($j\in\Z_{\ge 0}$) we get
\begin{equation}\label{eq:casimirexp}
C_{L(\fm+\un{j})}=\sum_{k=0}^{j-1}\xi^k\car(\fm+\un{k+1}).
\end{equation}
We can assume $j_1<j_2$. By assumption we have
\begin{align*}
0&=\sum_{k=0}^{j_2-1}\xi^k\car(\fm+\un{k+1})-
\sum_{k=0}^{j_1-1}\xi^k\car(\fm+\un{k+1})
=\sum_{k=j_1}^{j_2-1}\xi^k\car(\fm+\un{k+1})=\\
&=\xi^{j_1}\sum_{k=0}^{j_2-j_1-1}\xi^k\car(\fm+\un{j_2}-\un{(j_2-j_1)}+\un{k+1}).
\end{align*}
By Proposition \ref{prop:dimLm} this means that $\dim L(\fm+\un{j_2})\le j_2-j_1$.
But this contradicts Theorem \ref{thm:mainresult_ny} which says that
$\dim L(\fm+\un{j_2})=2j_2+1$.
\end{proof}

\begin{thm}\label{thm:ssodd}
Let $V$ be a finite-dimensional weight module over $A$
with support in $G_0+\un{\Z}$. Then $V$ is semisimple.
\end{thm}
\begin{proof}
We follow the idea of the proof of Proposition 12 in \cite{KliSch}, Chapter 3.
Writing
\[V=\oplus_{\fm\in G_0}\big(\oplus_{j\in\Z}V_{\fm+\un{j}}\big)\]
and noting that $\oplus_{j\in\Z}V_{\fm+\un{j}}$ are submodules,
we can reduce to the case when $\supp(V)$ is contained in $\fm+\Z$
for a fixed $\fm\in G_0$.

Let $\la_1,\ldots,\la_k$ be the generalized eigenvalues of the
Casimir operator $C_V$, i.e. the elements of the set
\[\{\la\in\K\, |\, \ker(C_V-\la\id)^p\neq 0\text{ for some $p>0$}\}.\]
Then each generalized eigenspace $\sum_p \ker(C_V-\la_i\id)^p$ is
invariant under $A$, hence they are submodules. It suffices to
prove that each such submodule is semisimple. Let $V$ be one of them.
Let $V_1=\{v\in V \, |\, X_+v=0\}$. Then $V_1$ is invariant under $R$
and since $V$ is a weight module, $V_1=\oplus_{\fn\in G}(V_1\cap
V_\fn)$. Now if $0\neq v\in V_1\cap V_\fn$, then $v$ is a highest
weight vector of $V$ and generates a submodule isomorphic to
$L(\fn)$. Hence if $V_1\cap V_\fn\neq 0$ for more than one $\fn\in
G$, $C_V$ will have two different eigenvalues by Proposition
\ref{prop:separation} which is impossible. Here we used that the
restriction of $C_V$ to a submodule $W$ coincides with $C_W$.
Hence $V_1$ is contained in a single weight space, say $V_{\fn}$.
Let $v_1,\ldots,v_k$ be a basis for $V_1$. Then each $v_i$
generates a simple submodule isomorphic to $L(\fn)$. We will
show that the sum of these submodules is direct. Vectors of
different weights are linearly independent so it suffices to show
that if
\[\sum_{i=1}^k \la_i (X_-)^mv_k=0\]
then all $\la_i=0$. Assume the sum was nonzero and act by $X_+$ $m$ times.
In each step we get a nonzero result because we have not reached the
highest weight $\fn$ yet. But then, using \eqref{eq:xnxn}, we have a
linear relation
among the $v_k$ -- a contradiction. We have shown that $V$ contains
the direct sum $V'$ of $k$ copies of $L(\fn)$. Now $X_+$ acts injectively
on $V/V'$. This is only possible in a torsion-free finite-dimensional
weight $A$-module if it is $0$-dimensional.
Thus $V$ is semisimple.
\end{proof}

We now turn to the general case.
Assume now that $A$ has an even-dimensional irreducible representation.
By Lemma \ref{lem:mainresult}b), $G_{1/2}\neq\emptyset$. We 
fix $\fm_{1/2}\in G$.
Then $G_{1/2}=G_0+\fm_{1/2}$ by Proposition \ref{prop:G12coset}.
\begin{thm}
Any finite-dimensional weight module $V$ over $A$ is semisimple.
\end{thm}
\begin{proof}
By Corollary \ref{cor:weightsupp} and Theorem \ref{thm:mainresult_ny},
\[\supp(V)\subseteq (G_0+\un{\Z})\cup(G_{1/2}+\un{\Z})\]
Thus we have a decomposition
\[V=\Big(\bigoplus_{\fm\in G_0} V_{\fm+\un{\Z}}\Big)
\oplus\Big(\bigoplus_{\fm\in G_0} V_{\fm+\fm_{1/2}+\un{\Z}}
\Big)\]
where $V_{\fn+\Z}:=\oplus_{j\in\Z}V_{\fn+\un{j}}$ for $\fn\in G$
are submodules.
It remains to prove that a weight module $V$ with support
in $\fm+\fm_{1/2}+\Z$ is semisimple.
By Lemma \ref{lem:tensorweight},
\[\supp\big(V\otimes L(\fm_{1/2})\big)\subseteq \fm+\fm_{1/2}+\fm_{1/2}+\Z=
\fm'+\Z\]
where $\fm':=\fm+\fm_{1/2}+\fm_{1/2}-1\in G_0$ by
Lemma \ref{lem:dim2special}b). Hence $V\otimes L(\fm_{1/2})$
is semisimple by Theorem \ref{thm:ssodd}. By the Clebsch-Gordan formula
\eqref{eq:clebschgordan}, the tensor product of two
semisimple modules is semisimple again. Therefore
$V\otimes L(\fm_{1/2})\otimes L(\un{1}-\fm_{1/2})$ is semisimple,
where $\dim L(\un{1}-\fm_{1/2})=2$ by Lemma \ref{lem:dim2special}a).
On the other hand, by \eqref{eq:clebschgordan} again we have
\[V\otimes L(\fm_{1/2})\otimes L(\un{1}-\fm_{1/2})\simeq
V\otimes\big( L(0)\oplus L(\fm)\big)\simeq
\big(V\otimes L(0)\big)\oplus\big( V\otimes L(\fm)\big).\]
Finally, it is easy to verify the isomorphism
 $V\simeq V\otimes L(0)$, $v\mapsto v\otimes e$
where $0\neq e\in L(0)$ is fixed. Thus $V$ is isomorphic to a submodule
of the semisimple module $V\otimes L(\fm_{1/2})\otimes L(\un{1}-\fm_{1/2})$
and is therefore itself semisimple.
\end{proof}

\end{document}